# Effective Implementation of GPU-based Revised Simplex algorithm applying new memory management and cycle avoidance strategies

Arash Raeisi Gahrouei, Mehdi Ghatee[1]

Department of Computer Sciences, Amirkabir University of Technology, Tehran, Iran

**Abstract**

Graphics Processing Units (GPUs) with high computational capabilities used as modern parallel platforms to deal with complex computational problems. We use this platform to solve large-scale linear programing problems by revised simplex algorithm. To implement this algorithm, we propose some new memory management strategies. In addition, to avoid cycling because of degeneracy conditions, we use a tabu rule for entering variable selection in the revised simplex algorithm. To evaluate this algorithm, we consider two sets of benchmark problems and compare the speedup factors for these problems. The comparisons demonstrate that the proposed method is highly effective and solve the problems with the maximum speedup factors 165.2 and 65.46 with respect to the sequential version and Matlab Linprog solver respectively.

**Keywords:** GPU-computing; Revised simplex algorithm; Degeneracy conditions; Multiple GPUs; Large scale problems;

---

[1] Corresponding author

Department of Computer Sciences, Amirkabir University of Technology, No 424, Hafez Ave, Tehran 15875-4413, Iran Fax: +98216497930, Email: ghatee@aut.ac.ir, URL: http://www.aut.ac.ir/ghatee

## 1: Introduction

To solve linear programing problems, the modern computers and efficient algorithms with reasonable computational time have been pursued [1]. The most important groups of solution methods for these problems are Simplex algorithm [2], Interior-point algorithm [3], Khachiyan algorithm [4] and so on. Simplex algorithm is one of the best solutions for real problems, although it can fail in degeneracy and cycling. One of the most popular versions of Simplex algorithm is Revised Simplex algorithm [1].

It should be noted that different applications are developed for solving linear programing problems such as Groubi, AMPL, CPLEX and Matlab see e.g. [5-7]. Nowadays, graphics processing units (GPUs) along with their distinct memories can be used as coprocessors in modern devices for improving the speed and solving quality of the algorithms. Using such devices alongside CPUs can provide a high performance and efficiency in distributed systems [8]. Various frameworks are developed based on different programing languages for using these capabilities, for example, CUDA, a framework similar to C language, has been designed for NVIDIA GPUs [9]. Also this framework includes some libraries for different purposes, CULA, CUBLAS and NVIDIA LAPAC which are used in numerical calculations [10, 11].

Although parallel algorithm implementation on GPU needs to understand its architecture to adapt with algorithm, the primary issue is hierarchy memories management method in GPU [8, 12].

GPU has been used in many computational projects in recent decades. Bolz et al. [13] used GPU to execute conjugate gradient method to solve spars matrix in 2003. In the same year, Krueger and Westermann solved linear systems by GPU instead of CPU [14]. In [15] researchers have tried to exploit the capability of GPU architecture to solve difficult problems in operating research area. Among computational problems, GPU is widely used to implement optimization algorithms. In Table 1, we have presented the most outstanding papers in this category.

Table 1. The most important references that have used GPU to implement optimization algorithms.

| Reference | Type of Optimization Problem | Approach | Speedup |
|---|---|---|---|
| Jung and O'Leary [16] | Linear Optimization Problem | interior point method | Do not report |
| Bieling et al. [17] | Linear Optimization Problem | revised simplex algorithm | 18 |
| Spampinato and Elstery [18] | Linear Optimization Problem | revised simplex algorithm | 2.5 |
| Lalami et al. [19] | Linear Optimization Problem | simplex algorithm | 12.5 |
| Lalami et al. [20] | Linear Optimization Problem | revised simplex algorithm on two GPUs | 24.5 |
| [21-23] | combinatorial optimization problems | local search and bees swarm algorithms | 100-280 |
| Ploskas et al. [24] | Linear Optimization Problem | considered scaling technique for linear programing problems | 7 |
| Ploskas et al. [25] | Linear Optimization Problem | Using some pivoting rules for the revised simplex algorithm | 16 |
| Ploskas et al. [26] | Linear Optimization Problem | Using five different updating schemes for inverse computing in simplex algorithm | 5.5 |
| Ploskas et al. [27] | Linear Optimization Problem | Simplex based Algorithm | 58 |
| Chen et al. [28] | Linear Optimization Problem | Iterative linear solver | 40 |

One important concept in the improvement of optimization algorithms on GPU is the method of matrix multiplication. Hull et al. [29] devised a strategy to overcome bandwidth problem in matrix multiplication by using Texture memory. The matrix-matrix multiplication performance has been studied in [30]. Reference [31], shows usage of GPU for matrix-vector multiplication in a high performance way. In that study, the speedup was reported more than 32 times faster than sequential implementation. In [32], a high performance method was proposed to compute matrix-matrix multiplication.

To solve these problems, in continuation of these works, in this study, by using the GPUs capabilities, we tried to implement the Revised Simplex algorithm to solve large scale linear programing problems efficiently. For this aim, our contribution includes:

- Applying some strategies for GPU memory management
- Using tabu rule for entering variable to avoid cycling problems

To evaluate the proposed Simplex algorithm, various comparisons have been made on famous benchmark of linear programing problems [33] and random linear programing problems with different scales.

The outline of this study is as follows: in Section 2, the Revised Simplex algorithm and GPU architecture are discussed. In Section 3, we propose some strategies to implement the Revised Simplex algorithm on GPU. In Section 4, the experimental results are illustrated; Final section includes conclusion and future directions.

## 2: Preliminaries

### 2.1. Revised Simplex algorithm

Any linear programing problem can be replaced as the following problem:

$$Min\ z = C^T X$$

$$s.t.:$$

$$AX = b \quad b \geq 0$$

$$X \geq 0$$

$$C = (C', 0, \dots, 0) \in R^{n+m}, A = (A', I_m) \in R^{m \times (n+m)}$$

$$X = (X', x_{n+1}, \dots, x_{n+m})^T \in R^{(n+m)}$$

Where $X = (x_1, \dots, x_n)^T \in R^n$ is the decision variables, $C' = (C_1, C_2, \dots, C_n) \in R^n$ is the vector of variable costs and matrix $A' = \begin{pmatrix} a_{11} & \cdots & a_{1n} \\ \vdots & \ddots & \vdots \\ a_{m1} & \cdots & a_{mn} \end{pmatrix} \in R^{m \times n}$ is the constraints coefficients. $b = (b_1, \dots, b_m)^T \in R^m$ also denotes the right hand side of constraints. $m, n$ are the number of variables and the number of constraints.

To solve this problem, in 1974, Dantzig presented the Revised Simplex algorithm [2]. Revised Simplex algorithm contains the following steps: (see e.g. [1] for details),

Step 1: Decompose the matrix A as $A = (B, N)$, in which B is an invertible $m \times m$ matrix. With respect to this decomposition, consider X as $X = (X_B, X_N)^T$ in which $X_B$ is an m vector. Similarly, we can decompose $C = (C_B, C_N)$. In this case, we can rewrite the standard linear program as the following.

$$Min z = C_B^T X_B + C_N^T X_N$$

$$s.t.:$$

$$BX_B + NX_N = b$$

$$X_B, X_N \geq 0$$

B is a basic matrix, N is a non-basic matrix, $X_B$ is vector of basic variables and $X_N$ is vector of non-basic variables.

Step 2: If $X_B = B^{-1}b = \bar{b}$, $X_B \geq 0$ and $X_N = 0$, then $X = (X_B, X_N)$ is a feasible basic solution for the linear programing problem with the objective value $z = C_B B^{-1} \bar{b}$.

Step 3: Calculate the reduced costs for each non-basic variable as the following:

$$z_j - c_j = C_B B^{-1} a_j - c_j = W a_j - c_j$$

where $a_j$ is a column of matrix N with respect to $j^{th}$ non-basic variable.

Step 4: Determine the maximum of reduced costs, $z_k - c_k = Max\ (z_j - c_j)$. If this value is not positive, the basic solution $X = (X_B, X_N)$ is optimal and the process terminates.

Step 5: Calculate $Y_k = B^{-1} a_k$, if $Y_k \leq 0$, then the objective function is unbounded and the process terminates.

Step 6: Denote the r$^{th}$ variable of the basic solution with $X_{B_r}$ where the index 'r' is determined from the following test:

$$\frac{\overline{b_r}}{y_{rk}} = \min\left\{\frac{\overline{b_i}}{y_{ik}}\ \bigg|\ i = 1, \ldots, n, y_{ik} > 0\right\}$$

Step 7: The k$^{th}$ non-basic variable $X_{N_k}$ enters the basic variables instead of $X_{B_r}$.

Step 8: Update the parameters and go to Step 2. For this aim, it is possible to use the tableau given in Figure 1. Then, r$^{th}$ row is divided by $y_{rk}$. For other rows such as 'i', the modified r$^{th}$ row is multiplied by $-y_{ik}$ and the results are added to row 'i'.

| W | $C_B \overline{b}$ | $z_k - c_k$ |
|---|---|---|
| $B^{-1}$ | $B^{-1} b$ | $Y_k$ |

Figure 1. The Revised Simplex algorithm tableau [1].

The Revised Simplex algorithm gets basic matrix, non-basic matrix, coefficient matrix, and right hand side matrix and then, it sends the common value of objective function and common basic variables as desirable results. For further study, refer to [1].

### 2.2. Procedures for preventing cycling

Sometimes, the Revised Simplex algorithm drops in a degenerate extreme solution and it is very hard to leave this solution and move to another extreme solution. In the worst case, it is possible to return to a degenerate solution repeatedly and the algorithm cannot terminate which is called cycling. Although the Lexicographic rule is high computational strategy for cycling avoidance, it is the most important rule that is used by many researchers [1]. In the proposed version of our Revised Simplex algorithm, we use another approach namely Tabu rule. Based on the Tabu rule, if we do not have unitary leaving element in an iteration of the Revised Simplex algorithm, we check the improvement of the objective function for every leaving element and the leaving element is chosen based on the most noticeable improvement in the objective function. Then, the selected leaving variable is augmented to a Tabu list related to the entering variable. These leaving elements cannot be selected in the next iterations by the algorithm and this approach avoids cycling.

### 2.3. GPU architecture

GPU architecture belongs to single instruction, multiple data (SIMD) class of parallel computers. A GPU includes several processors known as Thread processor and each of them includes special registers. Several Thread processors are titled as SM that access a shared memory. Each GPU includes several types of memories and different kinds of access to these memories are possible by the CPU and the Thread processors.

Each running process on the GPU is called a Kernel which includes a Grid; each Grid includes several Blocks and each Block includes several Threads. The process is started on the CPU and then, it is transferred to the GPU by calling a Kernel. When all Threads finish their tasks in the Grid, the Kernel is terminated and the process will continue on the CPU again. By calling each Kernel, each SM should process a Block in 32 Threads group known as Warp; see [8, 12] for further information.

### 3. GPU-based Revised Simplex algorithm

In this section, the implementation of our proposed Revised Simplex algorithm is presented in details. For this aim, the applied techniques for memory allocation, finding the initial feasible solution and cycling prevention are discussed.

It should be noted that all of the steps of Revised Simplex algorithm can be implemented in parallel manner and so, using GPU is an acceptable framework for this implementation. We used CUBLAS algebraic library for vector scaling, matrix-vector multiplication and copying a vector into the deferent memories. The pivoting or matrix updating was programmed completely by CUDA. The memory management method is an important issue in GPU to improve the performance of software. There are different implementation procedures based on the place which the whole tableau or a part of it, is kept on the Global memory of the GPU which are used. The flowchart of Figure 2 is implemented in which the free accessible Global memory of GPU is queried. Based on these amounts, the following cases are taken into consideration:

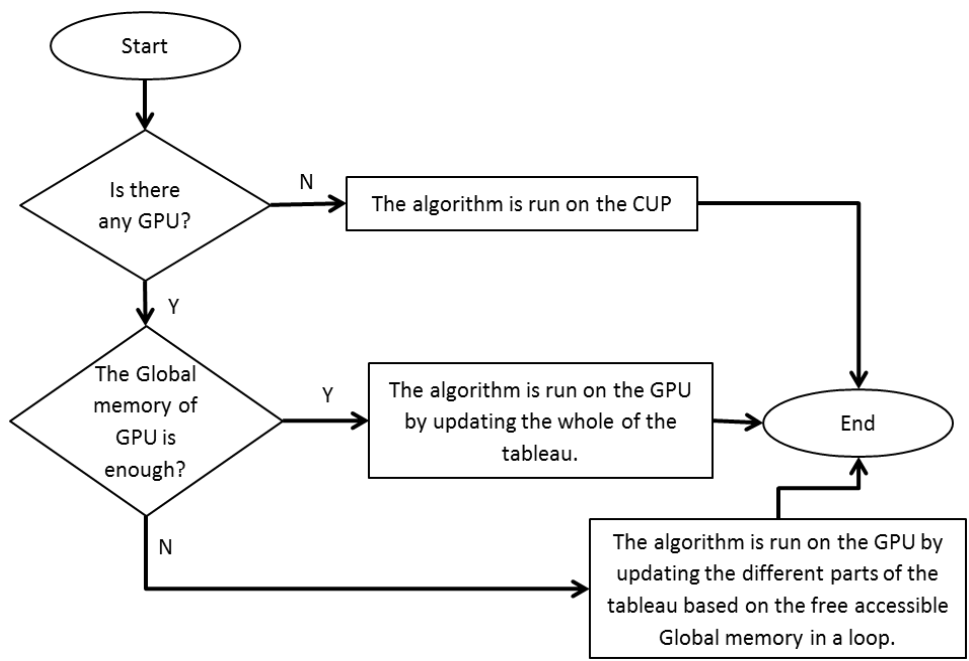

Figure 2. The flowchart of available Global memory on the GPU.

**Case 1: The GPU's Global memory is sufficient for keeping the whole Simplex tableau**

In this case, although the GPU Global memory size is sufficient to save all the tableau's elements, the tableau is transferred to GPU at once and it is kept in the GPU's Global memory up to algorithm termination. In this case, computational overhead of transferring between the main memory and GPU's memory is minimal. The details of implementation of the algorithm in this case are presented in Figure 3.

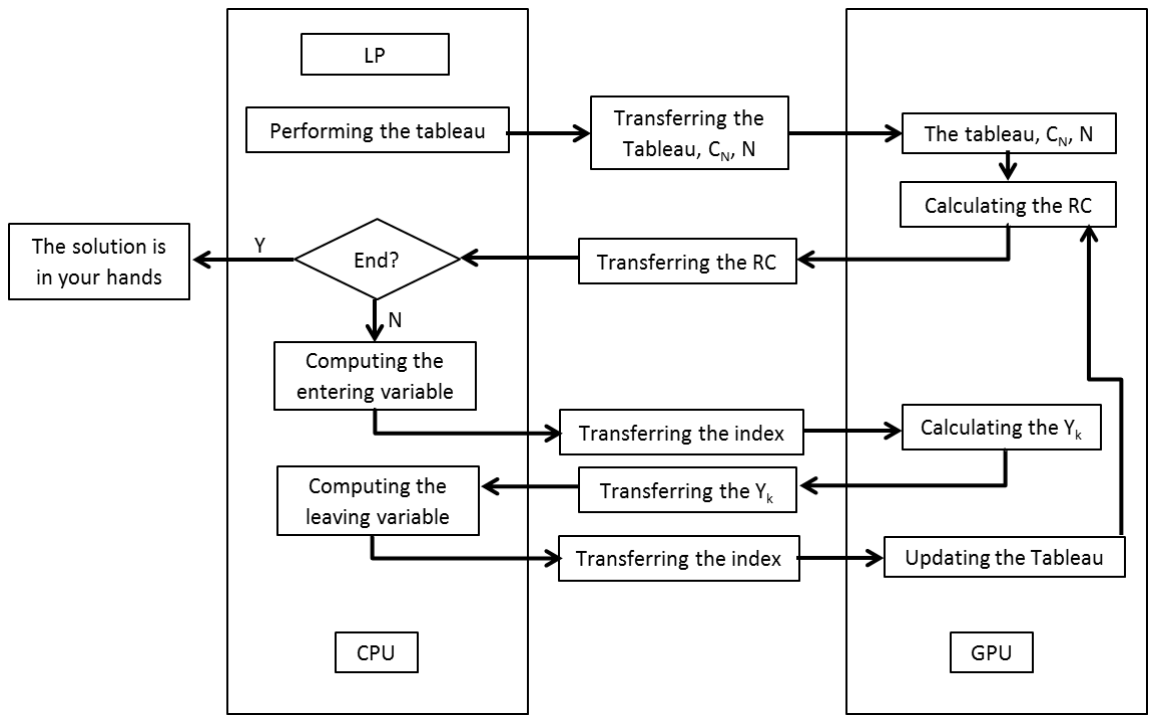

Figure 3. The implementation of algorithm on the GPU with sufficient Global memory.

**Case 2: The GPU's Global memory is insufficient for keeping the whole Simplex tableau**

In this case, we should unfortunately update the tableau through several steps according to the free accessible GPU's Global memory size. In this case, the tableau is divided into several parts according to the free accessible Global memory of the GPU and each part is updated separately. Therefore after updating a part of the tableau, it is transferred to the system's main memory and another part is sent to the GPU for updating. Thus in this case the performance is low especially when the number of transferring between the main memory and the GPU memory increases. To decrease the number of these transfers, the last part of the tableau is kept for the next updating. In addition to accelerating the process, the parts of the tableau are chosen in a way that several consecutive rows are considered in a part every time. In this case, r[th] row of tableau is transferred one time; see Figure 4 for details.

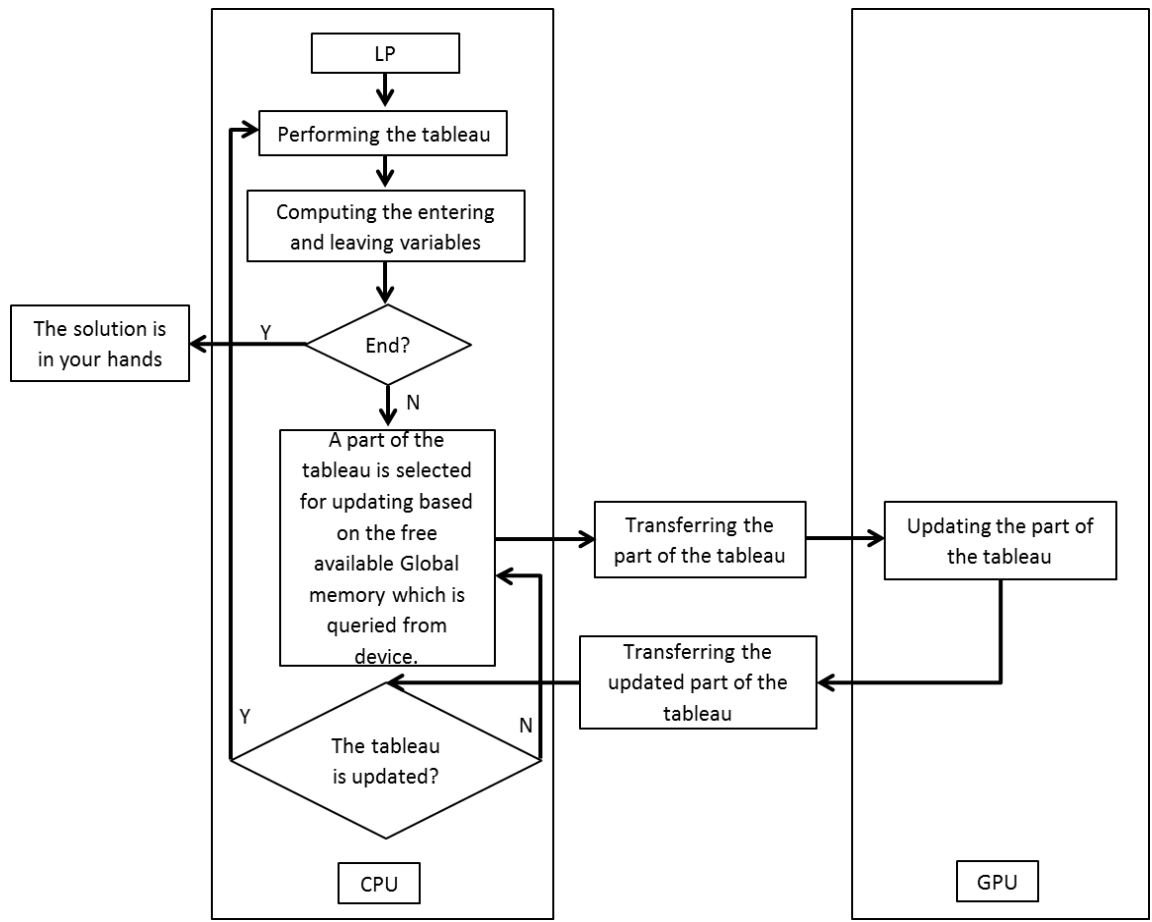

Figure 4. The implementation of algorithm on the GPU with insufficient Global memory.

## 3.1. Tableau updating

To update the Revised Simplex tableau in GPU, after identifying the entering variable 'k' and leaving row 'r', in the first step r[th] row elements of tableau are divided by pivot element ($y_{rk}$), then r[th] row is used to update other rows. It is necessary to avoid conditional command 'if i==r' in updating process inside the Kernel. For this aim, at first, r[th] row is divided by $y_{rk}$, Then, r[th] element of $Y_k$ is replaced with zero. Then, for i=1,...,m, the following steps are taken into account:

- Multiply row 'r' by $-y_{ik}$
- Add the results with row 'i'

After all of the updating steps, again r$^{th}$ element of $Y_k$ is replaced with 1. It is a simple exercise to show that the result of this version of updating is equivalent with updating process presented in Step 8 of Revised Simplex algorithm.

The summarization of this updating procedure is illustrated in Figure 5.

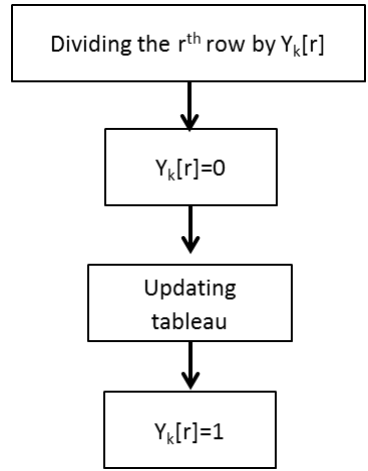

Figure 5. The summarization of updating process.

### 3.2. Allocation memory

The most important part of proposed algorithm is updating process. If 'kernel 1' is used for updating process and each block dimension is $16 \times 16$ which includes 256 Threads then $4 \times 256$ Global memory accesses are performed for each block updating. However, if the parts of d_y and d_x arrays (arrays that present $Y_k$ and r$^{th}$ row of tableau respectively) which should be accessible frequently by Threads in each Block are saved into the Shared memory of GPU then $32 + 2 \times 256$ Global memory accesses are needed. With regard to efficiency of GPU's Shared memory, we used 'kernel 2' for updating process.

**Kernel 1:**

```
__global__ void Updating_Kernel (double* d_A, double* d_x, double* d_y, int m, int n){
    int i = blockIdx.y*blockDim.x + threadIdx.y;
    int j = blockIdx.x*blockDim.y + threadIdx.x;
    if (i < m && j < n)
    {       double temp = d_A[i*n + j];
            temp = ((-d_y[i])* d_x[j]) + temp;
            d_A[i*n + j] = temp;
    }
    __syncthreads();
}
```

**Kernel 2:**

```
__global__ void Updating_Kernel (double* d_A, double* d_x, double* d_y, int m, int n){
    int i = blockIdx.y*blockDim.x + threadIdx.y;
    int j = blockIdx.x*blockDim.y + threadIdx.x;
    int p = threadIdx.x;
    int q = threadIdx.y;
    __shared__ double d_xs[16];
    __shared__ double d_ys[16];
    if (p == q)
    {       if (i < m)
            {       d_ys[p] = d_y[i];
            }
            if (j < n)
            {       d_xs[p] = d_x[j];
            }
    }
    __syncthreads();
    double temp = -d_ys[threadIdx.y] * d_xs[threadIdx.x];
    if (temp != 0)
    {       if (i < m && j < n)
            {       d_A[i*n + j] = temp + d_A[i*n + j];
            }
    }
    __syncthreads();
}
```

### 3.3. Block size determination

There are some restrictions on the Shared memories, the registers memories, the number of Threads in a Block, the number of Blocks in a Grid, and the maximum number of allocated Blocks to each SM on each GPU. So, block size determination is necessary for implementation of GPU based algorithm before running each Kernel. As a result, the high performance can be provided by choosing desirable block size.

To maximize GPU occupancy, in our implementation, each Block can be involved in 256 Threads according to CUDA occupancy calculator [9], so the blocks size is chosen 16 × 16 which is included 256 Threads.

**3.4. Cycling prevention**

As mentioned in Section 2, we can use Lexicographic rule to implement the Revised Simplex algorithm to prevent cycling. In the proposed version of our Revised Simplex algorithm, we use another approach namely Tabu rule. Based on the Tabu rule, if we do not have unitary leaving element in an iteration of the Revised Simplex algorithm, we check the improvement of the objective function for every leaving elements and the leaving element is chosen based on the most noticeable improvement in the objective function. Then, the selected leaving variable is augmented to a Tabu list related to the entering variable. These leaving elements cannot be selected in the next iterations by the algorithm and this approach avoids cycling.

**3.4. Finding the initial feasible solution**

Recall that the Revised Simplex method starts with a basic feasible solution and moves to an improved basic feasible solution, until the solution is reached. However, in order to initialize the Revised Simplex method, an initial feasible solution must be available. To find such a solution, we defined artificial variables and used the two-phase method to eliminate them in our algorithm. For further study, refer to [1].

**4. Computational results**

The computational comparison of GPU-based implementation of Revised Simplex algorithm has been performed on an Intel core i7 CPU (4790K, 4.00 GHz) with 16 gigabyte of Main memory, running under windows 8.1, 64bit and the NVIDIA GeForce GTX Titan Black GPU with 6 gigabyte of Global memory that has 2880 cores with 980 MHz clock rate. The algorithm has been implemented using Microsoft Visual Studio 2013, CUDA V.6.5 and CUBLAS library. Execution time of the algorithm has been measured in seconds using clock function in Microsoft Visual Studio. The Simplex and

Interior-point algorithm of Linprog function in Matlab R2014b software are used for computational comparison.

In this study, 36 optimization linear problems form Netlib benchmark and 26 randomly generated feasible linear problems with different number of variables and constraints were used, out of randomly generated problems, 12 problems were of dense type (D). The remaining 14 problems were divided into two parts, in 7 of which (S20), 20 present of the elements in constraints coefficient matrix (A), and in the other 7 (S60), 60 present were considered zero. The specifications of the problems are given in Appendix A.

It is necessary to mention that dense linear problems have commonly happened in a lot of applications [34]. Also the reported results by most of GPU-based algorithms have been compared on dense linear programing problems. The optimum values of problems in this study checked by Linprog function of Matlab software. In order to prevent cycling, Tabu rule has been used throughout the implementations. Global memory along with Shared memory is used in order to update operation in the GPU-based implementation. In order to assess the efficiency of the proposed algorithm, each problem was solved using the proposed algorithm as well as the Simplex and Interior-point algorithms available in the Linprog function provided by Matlab software package. The runtime of each algorithm in achieving the optimal solution was calculated in seconds. Since the proposed and simplex algorithms include some iteration to update the Simplex tableau, time per iteration (TPI) was additionally obtained for these two algorithms. Given the parallel processing capability of Matlab, in order to demonstrate the efficiency of using GPU in the proposed method, the results of sequential version of the algorithm, in which merely the CPU is used as the processor, were considered in some of the instances. To analyse the performance of the proposed GPU-base algorithm, the following speedup factor has been considered:

$$S_p(n) = \frac{T^*(n)}{T_p(n)}$$

where n is the problem size, p is the number of processors, $T^*(n)$ is the execution time of the sequential or Matlab implementation and $T_p(n)$ is the execution time of the parallel implementation (GPU implementation). The obtained results for each instance of the studied set are presented in Tables 2 and 3.

As shown in Table 2, in the case of randomly generated spars linear problems, the Interior point algorithm of Matlab was unable to find a solution to the problem as the scale of the problem increased, and was ultimately terminated with an 'Out of Memory' error.

Table 2. The obtained results for randomly generated liner problems.

| Instance Name | Matlab - Simplex Total Run Time | Matlab-Simplex No. of Iterations | Matlab-Simplex Time Per Iteration | Matlab-Interior-Point Total Run Time | CPU Total Run Time | GPU Total Run Time | GPU No. of Iterations | GPU Time Per Iteration | Matlab-Simplex Total Runtime Speedup | Matlab-Simplex Time Per Iteration Speedup | Matlab-Interior-Point Total Run Tim Speedup | CPU Total Run Time Speedup |
|---|---|---|---|---|---|---|---|---|---|---|---|---|
| 200_370_D | 2.682403 | 2031 | 0.001321 | 0.882147 | 4.758 | 1.042 | 3243 | 0.000321 | 2.574283 | 4.110488 | 0.84659 | 4.566219 |
| 500_800_D | 12.74139 | 4104 | 0.003105 | 5.680032 | 134.659 | 5.457333 | 15797 | 0.000345 | 2.334728 | 8.986769 | 1.040807 | 24.67487 |
| 800_1700_D | 129.4957 | 16221 | 0.007983 | 30.87031 | 997.583 | 14.08333 | 37286 | 0.000378 | 9.194959 | 21.13577 | 2.191975 | 70.8343 |
| 1000_1600_D | 45.41527 | 5892 | 0.007708 | 52.45441 | 992.231 | 11.10467 | 25495 | 0.000436 | 4.089746 | 17.69655 | 4.723637 | 89.35261 |
| 1000_1800_D | 35.27596 | 4746 | 0.007433 | 38.00648 | 620.317 | 6.593333 | 14699 | 0.000449 | 5.350247 | 16.57043 | 5.76438 | 94.08246 |
| 1000_2300_D | 523.9995 | 37211 | 0.014082 | 64.98849 | 2792.397 | 40.542 | 61911 | 0.000655 | 12.92486 | 21.50414 | 1.602992 | 68.87665 |
| 1400_2000_D | 851.163 | 40698 | 0.020914 | 150.6467 | 6390.555 | 48.32333 | 87755 | 0.000551 | 17.61391 | 37.97998 | 3.117473 | 132.2457 |
| 1400_2300_D | 22.94765 | 1602 | 0.014324 | 120.0256 | 691.396 | 5.125333 | 8743 | 0.000586 | 4.4773 | 24.4351 | 23.4181 | 134.8978 |
| 1500_2700_D | 270.5778 | 15465 | 0.017496 | 203.3556 | 4634.803 | 31.09467 | 48463 | 0.000642 | 8.701744 | 27.26884 | 6.539888 | 149.0546 |
| 1800_3200_D | 1326.598 | 45517 | 0.029145 | 336.433 | 11496.75 | 69.59467 | 88328 | 0.000788 | 19.06178 | 36.99033 | 4.834177 | 165.1958 |
| 2500_4500_D | 1622.163 | 37690 | 0.04304 | 968.3331 | - | 103.012 | 81269 | 0.001268 | 15.74732 | 33.95514 | 9.400197 | - |
| 3500_5500_D | 17173.65 | 156438 | 0.109779 | 2325.749 | - | 679.664 | 305374 | 0.002226 | 25.26786 | 49.32399 | 3.42191 | - |
| 200_370_S20 | 3.397549 | 1447 | 0.002348 | 0.845835 | - | 0.888 | 2623 | 0.000339 | 3.826069 | 6.935576 | 0.952517 | - |
| 800_1700_S20 | 8.897195 | 1180 | 0.00754 | Out Of Memory | - | 1.084 | 2265 | 0.000479 | 8.207744 | 15.7547 | - | - |
| 1000_1800_S20 | 16.57024 | 1725 | 0.009606 | Out Of Memory | - | 2.812667 | 5847 | 0.000481 | 5.891292 | 19.96892 | - | - |
| 1400_2000_S20 | 577.7565 | 27811 | 0.020774 | Out Of Memory | - | 36.59733 | 66003 | 0.000554 | 15.78685 | 37.46644 | - | - |
| 1500_2700_S20 | 1138.792 | 41005 | 0.027772 | Out Of Memory | - | 57.63 | 90228 | 0.000639 | 19.7604 | 43.48108 | - | - |
| 2500_4500_S20 | 1382.592 | 28387 | 0.048705 | Out Of Memory | - | 123.0653 | 96132 | 0.00128 | 11.23462 | 38.04581 | - | - |
| 3500_5500_S20 | 9798.498 | 97075 | 0.100937 | Out Of Memory | - | 524.7987 | 230985 | 0.002272 | 18.67097 | 44.42661 | - | - |
| 200_370_S60 | 1.425057 | 951 | 0.001498 | 0.876771 | - | 0.765333 | 2546 | 0.000301 | 1.862008 | 4.984935 | 1.145607 | - |
| 800_1700_S60 | 7.88007 | 1180 | 0.006678 | 6.296796 | - | 7.298 | 17119 | 0.000426 | 1.079757 | 15.66472 | 0.862811 | - |
| 1000_1800_S60 | 4.800713 | 145 | 0.033108 | Out Of Memory | - | 0.756667 | 1124 | 0.000673 | 6.344555 | 49.18124 | - | - |
| 1400_2000_S60 | 79.4927 | 7577 | 0.010491 | Out Of Memory | - | 12.962 | 23227 | 0.000558 | 6.13275 | 18.79971 | - | - |
| 1500_2700_S60 | 325.7951 | 17026 | 0.019135 | Out Of Memory | - | 27.62133 | 42869 | 0.000644 | 11.79505 | 29.69824 | - | - |
| 2500_4500_S60 | 131.1377 | 4558 | 0.028771 | Out Of Memory | - | 22.43533 | 16467 | 0.001362 | 5.845144 | 21.11715 | - | - |
| 3500_5500_S60 | 962.0611 | 18648 | 0.051591 | Out Of Memory | - | 136.3113 | 60039 | 0.00227 | 7.057821 | 22.72332 | - | - |

Table 3. The average of execution time of Netlib benchmark linear problems in 3 runs.

| Instance Name | Matlab-Simplex Total Run Time | Matlab-Simplex No. of Iterations | Matlab-Simplex Time Per Iteration | Matlab-Interior-Point Total Run Time | GPU Total Run Time | GPU No. of Iterations | GPU Time Per Iteration | Matlab-Simplex Total Run Time Speedup | Matlab-Simplex Time Per Iteration Speedup | Matlab-Interior-Point Total Run Time Speedup |
|---|---|---|---|---|---|---|---|---|---|---|
| fit1d | 0.55982 | 1041 | 0.000537771 | 0.059628 | 0.47543735 | 2129 | 0.000223315 | 1.177484268 | 2.408130651 | 0.125417155 |
| recipe | 0.028999 | 32 | 0.000906219 | 0.017984 | 0.05995764 | 68 | 0.00088173 | 0.483658127 | 1.027773519 | 0.299945093 |
| e226 | 0.285409 | 485 | 0.000588472 | 0.042244 | 0.247986649 | 1025 | 0.000241938 | 1.150904701 | 2.432324369 | 0.17034788 |
| lotfi | 0.186245 | 365 | 0.00051026 | 0.035101 | 0.187316304 | 635 | 0.000294986 | 0.994280778 | 1.729776148 | 0.187388921 |
| grow7 | 0.183593 | 311 | 0.000590331 | 0.050574 | 0.153461539 | 485 | 0.000316416 | 1.196345363 | 1.865683283 | 0.329554887 |
| forplan | 0.258739 | 212 | 0.001220467 | 0.087983 | 0.241494808 | 383 | 0.000630535 | 1.071406055 | 1.935606222 | 0.364326673 |
| scsd1 | 0.127677 | 159 | 0.000803 | 0.022457 | 0.115564004 | 299 | 0.000386502 | 1.104816338 | 2.077610599 | 0.194325215 |
| sctap1 | 0.240532 | 87 | 0.002764736 | 0.027409 | 0.214931233 | 144 | 0.001492578 | 1.119111431 | 1.852322368 | 0.127524509 |
| fit2d | 19.505877 | 9462 | 0.002061496 | 0.426799 | 1.180433432 | 15425 | 7.65273E-05 | 16.52433459 | 26.93805338 | 0.361561261 |
| finnis | 0.389558 | 363 | 0.001073163 | 0.047035 | 0.157661735 | 542 | 0.000290889 | 2.470846834 | 3.6892534 | 0.298328569 |
| grow15 | 0.513084 | 735 | 0.000698073 | 0.091966 | 0.146437328 | 1280 | 0.000114404 | 3.50377876 | 6.101818794 | 0.628022931 |
| scsd6 | 0.275331 | 375 | 0.000734216 | 0.038092 | 0.144358444 | 727 | 0.000198567 | 1.907273264 | 3.697567101 | 0.263870952 |
| etamacro | 0.756408 | 529 | 0.001429883 | 0.116964 | 0.230862216 | 714 | 0.000323336 | 3.2764478 | 4.422275481 | 0.506639856 |
| pilot4 | 2.839675 | 2300 | 0.001234641 | 0.20485 | 0.196578898 | 5628 | 3.49287E-05 | 14.44547219 | 35.34744239 | 1.04207523 |
| grow22 | 1.099098 | 1284 | 0.000855995 | 0.140046 | 0.090765209 | 1853 | 4.89828E-05 | 12.10924332 | 17.47541112 | 1.542948027 |
| scrs8 | 0.768374 | 1044 | 0.00073599 | 0.071797 | 0.066198095 | 1413 | 4.68493E-05 | 11.60719209 | 15.70973412 | 1.084578045 |
| bnl1 | 6.898104 | 880 | 0.007838755 | 0.117783 | 0.817640872 | 1263 | 0.00064738 | 8.436593909 | 12.10842967 | 0.144052241 |
| ship04s | 0.385844 | 322 | 0.001198273 | 0.046383 | 0.134630654 | 805 | 0.000167243 | 2.86594464 | 7.164861601 | 0.344520351 |
| perold | 11.791097 | 4326 | 0.002725635 | 0.550619 | 0.494000214 | 6287 | 7.85749E-05 | 23.86860707 | 34.68838018 | 1.114612878 |
| ship04l | 0.567661 | 348 | 0.00163121 | 0.075884 | 0.10818541 | 510 | 0.000212128 | 5.24711236 | 7.689733631 | 0.701425453 |
| maros | 2.604774 | 1558 | 0.00167187 | 0.275002 | 0.184538191 | 3334 | 5.53504E-05 | 14.11509451 | 30.20521508 | 1.490217278 |
| shell | 1.664952 | 180 | 0.009249733 | 0.078167 | 0.772069605 | 394 | 0.001959568 | 2.156479142 | 4.720293233 | 0.101243462 |
| 25fv47 | 16.685754 | 5222 | 0.00319528 | 0.606141 | 0.591063082 | 8929 | 6.61959E-05 | 28.23007306 | 48.27007322 | 1.025509828 |
| fit1p | 4.100873 | 1233 | 0.003325931 | 0.14648 | 0.298366992 | 1686 | 0.000176967 | 13.74439234 | 18.79403527 | 0.490939025 |
| scsd8 | 2.405572 | 990 | 0.002429871 | 0.077413 | 0.188502509 | 1889 | 9.97896E-05 | 12.76148532 | 24.34994522 | 0.410673579 |
| sierra | 0.84316 | 391 | 0.002156419 | 0.153789 | 0.079287662 | 726 | 0.000109212 | 10.63418919 | 19.74532315 | 1.939633428 |
| pilot.we | 11.696551 | 7685 | 0.001521998 | 0.953649 | 0.379958275 | 10153 | 3.74233E-05 | 30.78377748 | 40.6698364 | 2.50987822 |
| pilotnov | 30.849332 | 4367 | 0.007064193 | 0.704805 | 1.117111752 | 9562 | 0.000116828 | 27.61526046 | 60.46648055 | 0.630917183 |
| sctap3 | 2.334103 | 747 | 0.003124636 | 0.143276 | 0.155215122 | 1822 | 8.51894E-05 | 15.03785827 | 36.67868509 | 0.923080165 |
| ship08l | 1.471129 | 409 | 0.003596892 | 0.168669 | 0.156284968 | 716 | 0.000218275 | 9.413118979 | 16.47871195 | 1.07924007 |
| czprob | 3.705756 | 1762 | 0.002103153 | 0.283288 | 0.165466977 | 2811 | 5.88641E-05 | 22.39574365 | 35.72896447 | 1.712051583 |
| d2q06c | 136.316149 | 53993 | 0.0025247 | 4.559453 | 2.082305108 | 84695 | 2.45859E-05 | 65.46406118 | 102.6888423 | 2.189618122 |
| 80bau3b | 9.932598 | 5400 | 0.00183937 | 0.492855 | 0.311043377 | 12007 | 2.59052E-05 | 31.93316027 | 71.00397321 | 1.584521764 |
| greenbea | 58.005099 | 14511 | 0.003997319 | 2.867399 | 1.512820369 | 21720 | 6.9651E-05 | 38.34235723 | 57.39066908 | 1.895399519 |
| maros-r7 | 83.668115 | 2551 | 0.032798163 | 25.32248 | 3.007690009 | 3412 | 0.000881504 | 27.81806461 | 37.20707035 | 8.419245309 |
| fit2p | 200.469882 | 13462 | 0.014891538 | 2.491421 | 4.159735612 | 31112 | 0.000133702 | 48.19293837 | 111.3785989 | 0.598937344 |

In Figure 6 diagram, the TPI of the proposed algorithm is compared to its sequential version (Mentioning that the number of iterations for GPU-based and CPU-based algorithms is same). As the graph shows clearly, there is a meaningful difference between TPI on CPU and GPU so that, TPI on GPU is at most 165 times less than the TPI on CPU. In addition, when the computational complexity of the problem increases, TPI on GPU increases gradually which shows that GPU-based implementation is more efficient than CPU-based implementation. Therefore, using GPU for solving linear programing problems with a great number of constraints and variables is more efficient than using CPU.

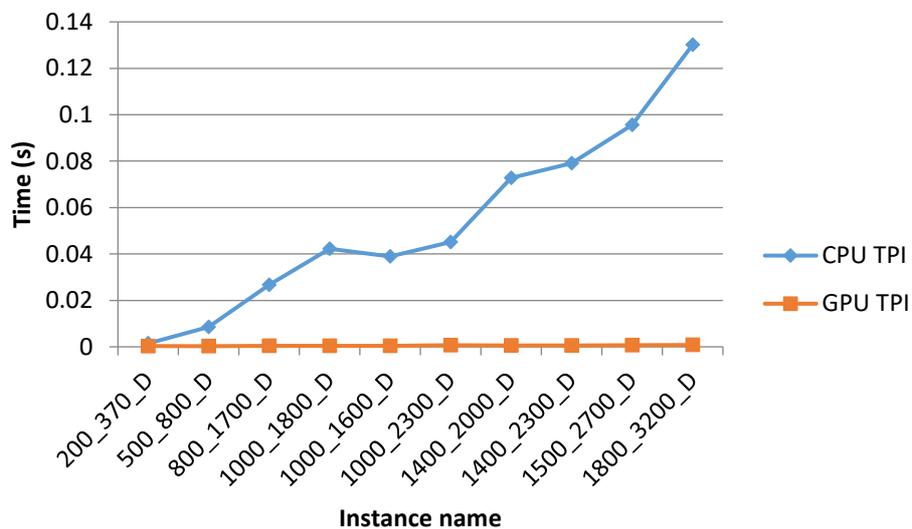

Figure 6. The comparison between TPI on CPU and TPI on GPU.

In Figure 7 diagram, the speedup of the proposed algorithm against the runtime to solve the problem using Matlab, is illustrated. As shown, as the solution runtime increases, the more appropriate speedups are expected, meaning that the more incapable Matlab in solving the problems, the more efficient the proposed algorithm in doing so. This indicates the increased efficiency of the algorithm in solving large-scale problems.

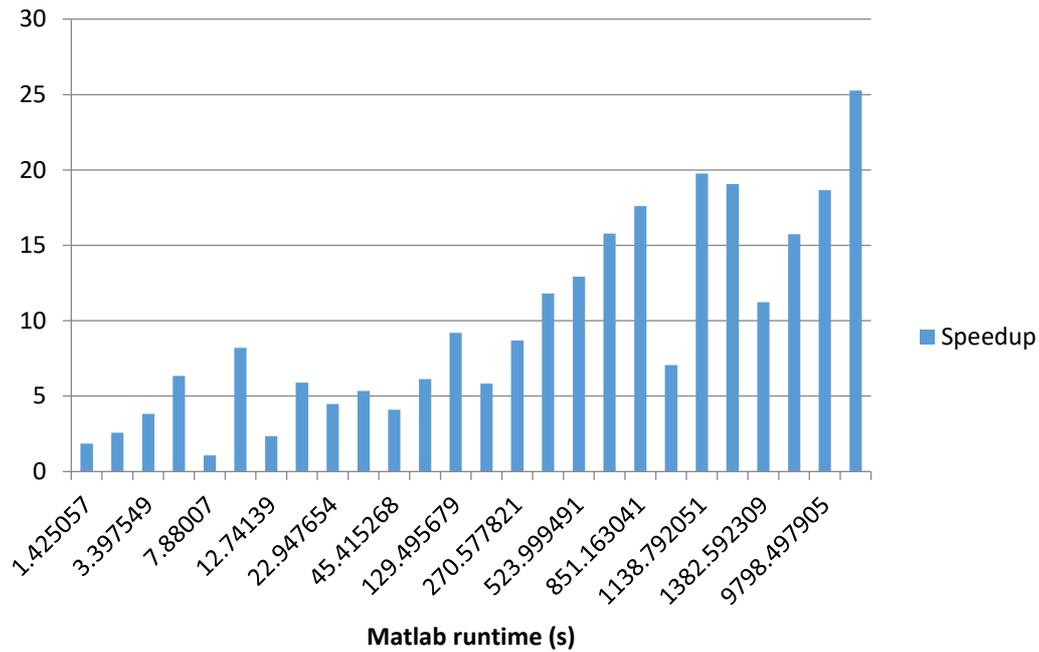

Figure 7. The speedup of the proposed algorithm against the runtime to solve the problem using Matlab.

In Figure 8 diagram, the maximum and mean of the obtained speedup per iteration for the three sets of randomly generated problems are compared to that of the Simplex algorithm of Matlab. As demonstrated, in case of sparser problems, the maximum and mean of speedups do not significantly vary. This is an indication of the efficiency of the algorithm in solving problems with different data structures.

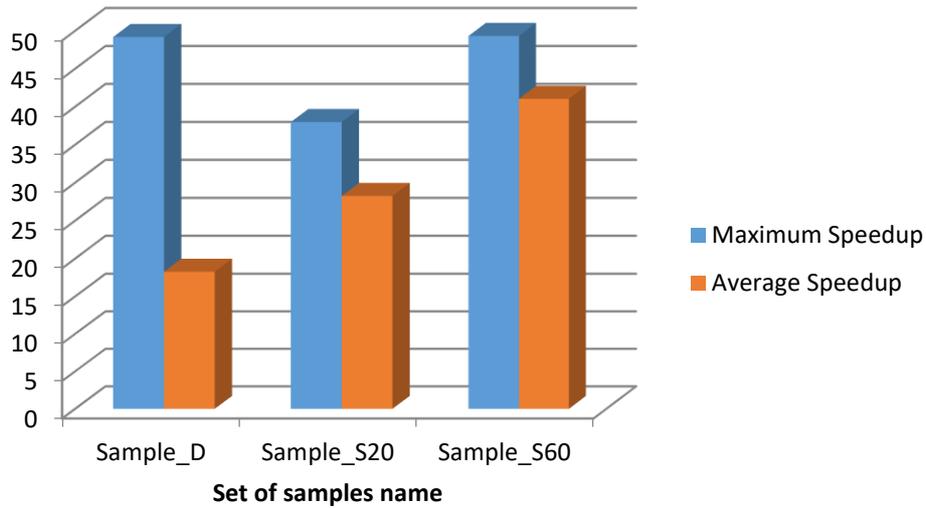

Figure 8. The maximum and mean of the obtained speedup per iteration for randomly generated problems.

In Figure 9 diagram, the obtained speedups per iteration for the Netlib benchmark problems are compared to that of Simplex algorithm provided by Matlab. According to the specifications of each problem presented in Appendix A, we can conclude that the proposed algorithm also outperforms the Simplex algorithm of Matlab in most of the problems with specific and regular numerical structures.

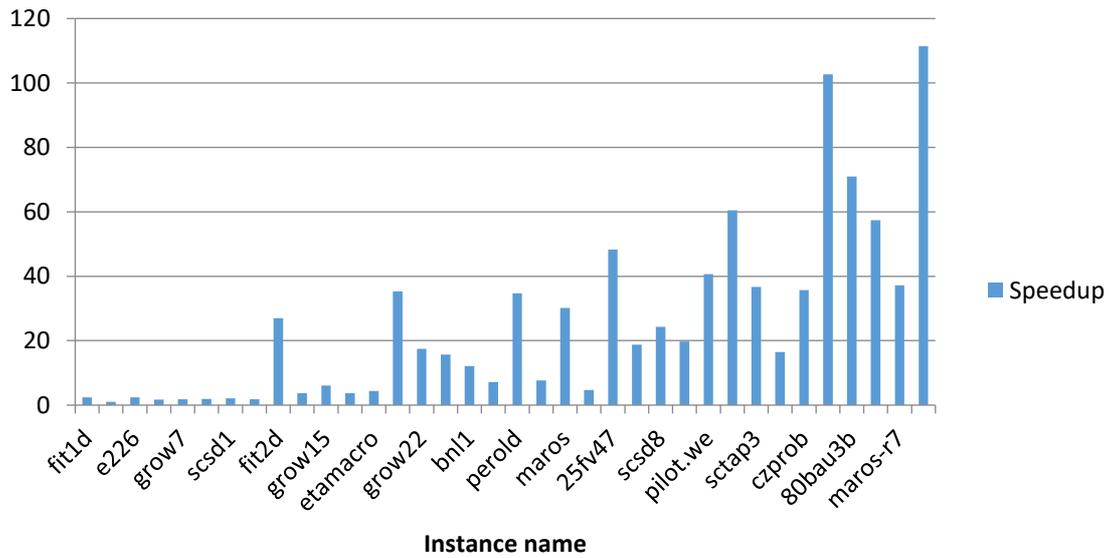

Figure 9. The comparison between Matlab simplex algorithm TPI speedup for the Netlib benchmark problems.

In Figure 10 diagram, the obtained speedups for Netlib benchmark problems are compared to that of Interior-point algorithm available in Matlab. In this case, as the computational complexity of the problem increases, the proposed algorithm better demonstrates its capability. It should be noted that due to the regular structure of such problems, each stage involved in the Interior-point algorithm is efficiently solved. However, in the case of randomly generated spars problems, which are irregular in structure, this algorithm proves inefficient so that it is not able to find a solution for many of the considered problems. This is demonstrated in Table 2.

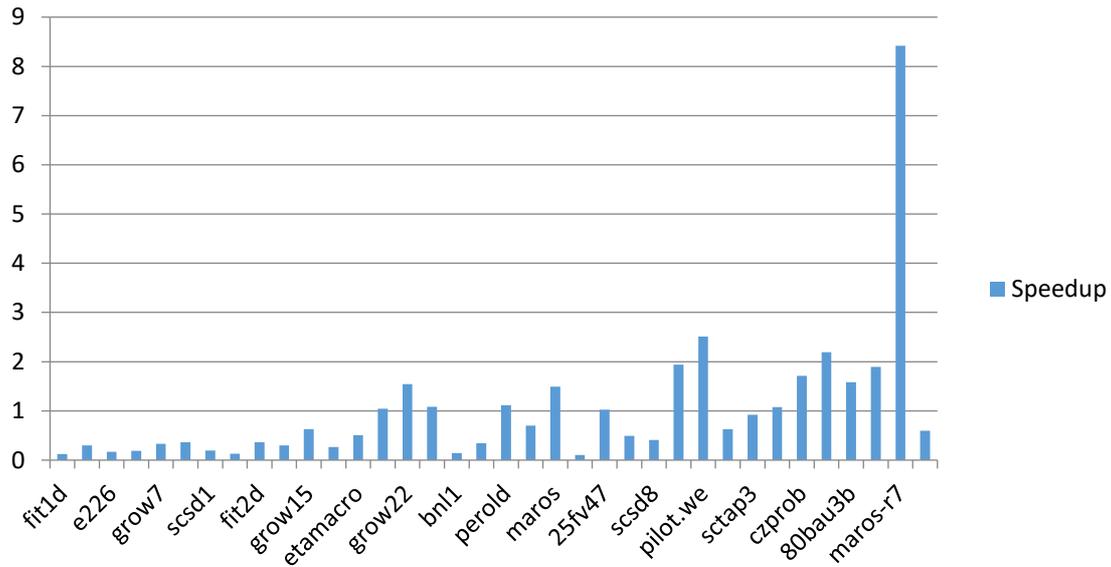

Figure 10. The comparison between Matlab interior-point algorithm TPI speedup for the Netlib benchmark problems.

In Table 4, the maximum and mean of achieved speedups are shown, according to which, it can be stated that compared to its sequential version, the proposed algorithm is highly capable in terms of computational power. Compared to the Simplex algorithm available in Matlab, the proposed algorithm, for both randomly generated problems as well as the Netlib benchmark problems, proves considerably efficient and finally it can easily compete with the Interior-point algorithm of Matlab, particularly in problems with irregular and sparse numerical structures.

Table 4. The maximum and mean of proposed algorithm speedup.

| Name of algorithm | Randomly generated problems | | Netlib benchmark problems | |
|---|---|---|---|---|
| | **Maximum speedup** | **Mean of speedup** | **Maximum speedup** | **Mean of speedup** |
| CPU-based version of proposed algorithm | 165.1958 | 93.3781 | - | - |
| Simplex algorithm of Matlab | 25.26785681 | 9.647452586 | 65.46406118 | 14.25552522 |
| Interior-point algorithm of Matlab | 23.4181 | 4.657538 | 8.419245309 | 1.023127 |

## 5. Conclusion

In this paper, with the help of CUDA, double precision implementation of the revised simplex algorithm is presented to solve the large scale linear programing problems on the GPU. For this aim, efficient memory management strategies have been considered. In addition to avoiding cycling in degenerate solutions, a new Tabu rule is considered. The speedup factor is considered to compare the results of this paper with the Matlab software. The value of the maximum speedup which has been achieved in this paper is 25.27 for randomly generated linear problem and 65.46 for Netlib benchmark linear problem. Thus, the presented implementation can be considered as the base structure in a series of problems which need to solve large-scale linear programing problems in real time. Our next works cover the application of such strategies for nonlinear programing problems.

## Appendix A

The specifications of the randomly generated liner problems.

| Instance Name | No. of Cons. | No. of Vars. | Optimum Value | No. of total constraints coefficients elements | No. of zeros constraints coefficients elements |
|---|---|---|---|---|---|
| 200_370_D | 200 | 370 | 195.7626 | 74000 | 0 |
| 500_800_D | 500 | 800 | 497.7864 | 400000 | 0 |
| 800_1700_D | 800 | 1700 | 211.5725 | 1360000 | 0 |
| 1000_1600_D | 1000 | 1600 | 78.98614 | 1600000 | 0 |
| 1000_1800_D | 1000 | 1800 | 53.99017 | 1800000 | 0 |
| 1000_2300_D | 1000 | 2300 | 155.2717 | 2300000 | 0 |
| 1400_2000_D | 1400 | 2000 | 438.9385 | 2800000 | 0 |
| 1400_2300_D | 1400 | 2300 | 2.213694 | 3220000 | 0 |
| 1500_2700_D | 1500 | 2700 | 71.48206 | 4050000 | 0 |
| 1800_3200_D | 1800 | 3200 | 123.0974 | 5760000 | 0 |
| 2500_4500_D | 2500 | 4500 | 26.96455 | 11250000 | 0 |
| 3500_5500_D | 3500 | 5500 | 651.6318 | 19250000 | 0 |
| 200_370_S20 | 200 | 370 | 143.5523 | 74000 | 13313 |
| 800_1700_S20 | 800 | 1700 | 6.296796 | 1360000 | 245767 |
| 1000_1800_S20 | 1000 | 1800 | 12.35752 | 1800000 | 325839 |
| 1400_2000_S20 | 1400 | 2000 | 154.6458 | 2800000 | 507325 |
| 1500_2700_S20 | 1500 | 2700 | 184.7414 | 4050000 | 733334 |
| 2500_4500_S20 | 2500 | 4500 | 36.79041 | 11250000 | 2037390 |
| 3500_5500_S20 | 3500 | 5500 | 121.3917 | 19250000 | 3487473 |
| 200_370_S60 | 200 | 370 | 279.2534 | 74000 | 33049 |
| 800_1700_S60 | 800 | 1700 | 6.296796 | 1360000 | 612954 |
| 1000_1800_S60 | 1000 | 1800 | 0.672388 | 1800000 | 810309 |
| 1400_2000_S60 | 1400 | 2000 | 84.73411 | 2800000 | 1261915 |
| 1500_2700_S60 | 1500 | 2700 | 114.8599 | 4050000 | 1824879 |
| 2500_4500_S60 | 2500 | 4500 | 9.336604 | 11250000 | 5072973 |
| 3500_5500_S60 | 3500 | 5500 | 15.3367 | 19250000 | 8679441 |

The specifications of the Netlib benchmark liner problems.

| Instance Name | No. of Cons. | No. of Vars. | Optimum Value | No. of total constraints coefficients elements | No. of zeros constraints coefficients elements |
|---|---|---|---|---|---|
| fit1d | 25 | 1026 | -9146.3781 | 25650 | 11211 |
| recipe | 158 | 180 | -266.616 | 28440 | 27272 |
| e226 | 256 | 282 | -18.751929 | 72192 | 68430 |
| lotfi | 248 | 308 | -25.264706 | 76384 | 74121 |
| grow7 | 280 | 301 | -47787812 | 84280 | 78778 |
| forplan | 252 | 421 | -664.21896 | 106092 | 97503 |
| scsd1 | 154 | 760 | 8.666667 | 117040 | 112114 |
| sctap1 | 420 | 480 | 1412.25 | 201600 | 199130 |
| fit2d | 26 | 10500 | -68464.293 | 273000 | 133470 |
| finnis | 544 | 614 | 172791.066 | 334016 | 331031 |
| grow15 | 600 | 645 | -106870941 | 387000 | 375162 |
| scsd6 | 294 | 1350 | 50.5 | 396900 | 387978 |
| etamacro | 672 | 688 | -755.71523 | 462336 | 457883 |
| pilot4 | 697 | 1000 | -2581.1376 | 697000 | 688590 |
| grow22 | 880 | 946 | -160834336 | 832480 | 815098 |
| scrs8 | 874 | 1169 | 904.296954 | 1021706 | 1015080 |
| bnl1 | 875 | 1175 | 1977.62956 | 1028125 | 1019564 |
| ship04s | 756 | 1458 | 1798714.7 | 1102248 | 1094303 |
| perold | 1120 | 1376 | -9380.7553 | 1541120 | 1529600 |
| ship04l | 756 | 2118 | 1793324.54 | 1601208 | 1589963 |
| maros | 1169 | 1443 | -58063.744 | 1686867 | 1672477 |
| shell | 1070 | 1775 | 1208825346 | 1899250 | 1891078 |
| 25fv47 | 1337 | 1571 | 5501.84589 | 2100427 | 2082793 |
| fit1p | 1254 | 1677 | 9146.37809 | 2102958 | 2081970 |
| scsd8 | 794 | 2750 | 905 | 2183500 | 2165542 |
| sierra | 1755 | 2036 | 15394362.2 | 3573180 | 3560156 |
| pilot.we | 1305 | 2789 | -2720107.5 | 3639645 | 3621361 |
| pilotnov | 1676 | 2172 | -4497.2762 | 3640272 | 3615331 |
| sctap3 | 2100 | 2480 | 1424 | 5208000 | 5195168 |
| ship08l | 1476 | 4283 | 1909055.21 | 6321708 | 6299020 |
| czprob | 1819 | 3523 | 2185196.7 | 6408337 | 6388830 |
| d2q06c | 3678 | 5167 | 122784.211 | 19004226 | 18947311 |
| 80bau3b | 2262 | 9799 | 987224.192 | 22165338 | 22142075 |
| greenbea | 4591 | 5405 | -72555248 | 24814355 | 24756292 |
| maros-r7 | 6272 | 9408 | 1497185.17 | 59006976 | 58711010 |
| fit2p | 6000 | 13525 | 68464.2933 | 81150000 | 81043434 |